\documentclass[12pt,a4paper]{amsart}

\usepackage[pdftex, pdfstartview=FitH,colorlinks=true,linkcolor=blue,%
urlcolor=blue,citecolor=blue,pagebackref=true]{hyperref}

\usepackage{amsmath}
\usepackage[abbrev]{amsrefs}

\usepackage{xcomment}

\newtheorem {Theorem}   {Theorem}
\newtheorem {thm}    [Theorem]{Theorem}

\newtheorem {Conjecture}[Theorem]{Conjecture}

\newtheorem {Definition} [Theorem]    {Definition}

\newtheorem {Question} [Theorem]    {Question}
\newcounter{AbcT}

\theoremstyle{definition}

\renewcommand{\a}{\alpha}
\renewcommand{\b}{\beta}
\renewcommand{\d}{\delta}

\newcommand{\e}{\varepsilon}
\newcommand{\f}{\varphi}

\renewcommand{\l}{\lambda}
\renewcommand{\o}{\omega}

\newcommand{\s}{\sigma}

\newcommand{\R}{{\mathbf R}}

\newcommand{\Z}{{\mathbf Z}}

\newcommand{\wt}{\widetilde}

\newcommand{\be}{\begin{equation}}
\newcommand{\ee}{\end{equation}}
\newcommand{\bea}{\begin{eqnarray}}
\newcommand{\eea}{\end{eqnarray}}
\newcommand{\bean}{\begin{eqnarray*}}
\newcommand{\eean}{\end{eqnarray*}}

\DeclareMathOperator{\supp}{supp}

\DeclareMathOperator{\Sim}{Sim}

\DeclareMathOperator{\dist}{dist}

\title[Entropy rates]{Entropy rates in the dimension theory of self-similar measures}
\author{P\'eter P. Varj\'u}

\begin{document}

\begin{abstract}
This is a survey paper on the dimension theory of self-similar
measures on the real line focusing on the role of entropy rates.
\end{abstract}

\maketitle

\section{Introduction}\label{sc:intro}

A (self-similar) iterated function system on the real line, IFS for short, is a finite collection
\[
\Phi=\{\f_i:i\in\Lambda\}
\]
of contractive similarities of $\R$.
A contractive similarity is a map $x\mapsto \l x +\tau$, where $\l\in(-1,1)\backslash\{0\}$,
$\tau\in\R$.
We call $\l$ the scaling factor and $\tau$ the translation parameter of the similarity.
Given such an IFS, and a probability vector $\{p_i:i\in\Lambda\}$, there
is a unique self-similar measure, that is a probability measure $\mu$ on $\R$ such that
\[
\mu=\sum_{i\in\Lambda} p_i \f_i(\mu).
\]
Here $\f_i(\mu)$ denotes the push-forward of $\mu$ under $\f_i$.
In other words, $\mu$ is the unique stationary measure for the Markov chain on $\R$ with
the transition $\f_i$ executed with probability $p_i$.

A measure $\mu$ on $\R$ is called exact dimensional, if there is a number $\a\in[0,1]$
such that for $\mu$-almost every $x$, we have $\mu([x-r,x+r])=r^{\a+o_{r\to 0}(1)}$.
If this is the case, we write $\dim \mu=\a$.
It is known (see Feng and Hu \cite{FH}) that self-similar measures are exact dimensional.

There are two natural upper bounds for the dimension of a self-similar measure on $\R$.
One of them is $1$, the dimension of the ambient space.
The other is a quantity depending on $\Phi$ and $p$, called the similarity dimension
\begin{equation}\label{eq:sim_dim}
\frac{\sum_{i\in\Lambda} p_i\log p_i}{\sum_{i\in\Lambda}p_i\log|\l_i|},
\end{equation}
where $\l_i$ is the scaling factor of $\f_i$.
Observe that the numerator is the Shannon entropy of $p$ (times $-1$) and the denominator
is the average logarithmic contraction rate.
See \cite{BSS}*{Theorem 3.2.1} for a proof that this is an upper bound for $\dim \mu$.

When the IFS satisfies a suitable separation condition, for example, if $\f_i(\supp\mu)$
are pairwise disjoint for $i\in\Lambda$, then it is a classical result
that $\dim\mu$ equals the similarity dimension \eqref{eq:sim_dim}.
Under this condition, the similarity dimension is never greater than $1$, so the dimension
of the ambient space does not matter.
For these results, it is enough to assume the so-called open set condition, which is
a mildly weaker separation condition.
See \cite{BSS}*{Definition 1.5.2} for the precise definitions of the separation conditions,
and see \cite{BSS}*{Theorem 3.2.1} and the references therein for the result about
$\dim\mu$ under these conditions.

The ultimate aim of the dimension theory of self-similar measures is to understand to
what extent these separation conditions may be relaxed.
One constraint to this comes from the fact that it may be possible to realize the same measure
as the self-similar measure
associated to several different choices of an IFS and a probability vector, which 
may have different similarity dimensions.
A simple instance of this is that when $\Phi$ contains the same map twice, it may be
replaced by a single instance and the corresponding probabilities may be replaced by
their sum in $p$.
This process lowers the similarity dimension.

We may also replace $\Phi$ by the compositions $\f_i\circ\f_j$ of all pairs
of maps for $(i,j)\in\Lambda^2$ and $p$ by the products $p_ip_j$.
This way we get a new IFS and a new probability vector, which define
the same probability measure.
It is easy to see that the similarity dimension does not change when we perform 
this operation, but the new IFS may contain the same map more than once even
if the maps in the original IFS were distinct.
Of course, we may consider compositions of order higher than $2$, too.
This motivates the following definitions.

\begin{Definition}[Exact overlaps]\label{df:exact-overlaps}
An IFS $\Phi$ contains exact overlaps, if there is a number $k\in\Z_{\ge1}$
and $(i_1,\ldots,i_k)\neq(j_1,\ldots,j_k)\in\Lambda^k$
such that
\[
\f_{i_1}\circ\cdots\circ\f_{i_k}=\f_{j_1}\circ\cdots\circ\f_{j_k}.
\]
\end{Definition}

In other words, an IFS does not contain exact overlaps if and only if
the semigroup generated by $\Phi$ with the composition operation is free.

\begin{Definition}[Entropy rate]\label{df:entropy-rate}
Let $\Phi$ be an IFS and $p$ a probability vector.
Let $X_1,X_2,\ldots$ be a sequence of independent random elements of $\Phi$ distributed according to
$p$.
The entropy rate is defined as
\[
h(\Phi,p)=\lim_{N\to\infty}\frac{1}{N}H(X_1\circ\cdots\circ X_N),
\]
where $H(\cdot)$ denotes the Shannon entropy.
\end{Definition}

The limit always exists and equals to the infimum, because the sequence is subadditive.
This quantity is also known as the random walk entropy, or as Garsia entropy after Garsia \cite{Gar:entropy}, who first studied
it in the context of Bernoulli convolutions,
a special class of self-similar measures
to be discussed in Section~\ref{sc:Bernoulli}.

In the absence of exact overlaps, it is easy to see that the entropy rate $h(\Phi,p)$ equals
the Shannon entropy $H(p)$ of the probability vector $p$.
It is a good way to think about the gap $H(p)-h(\Phi,p)$ as a quantity measuring
the amount of exact overlaps the IFS contains.

If we replace $\Phi$ with the IFS containing all $N$-fold compositions, and then
collapse the maps that appear with higher multiplicity, the corresponding probability
vector will have Shannon entropy $H(X_1\circ\cdots\circ X_N)$, and the average
logarithmic contraction rate is multiplied by $N$, so the new similarity dimension is
\[
\frac{H(X_1\circ\cdots\circ X_N)}{N\sum_{i\in\Lambda}p_i\log|\lambda_i|^{-1}}.
\]
So the infimum of the similarity dimensions of the IFS and probability vector pairs
that we may obtain by taking compositions and collapsing maps is given by
\[
\frac{h(\Phi,p)}{\sum_{i\in\Lambda}p_i\log|\lambda_i|^{-1}}.
\] 

This motivates the following folklore conjecture.

\begin{Conjecture}\label{cn}
Let $\Phi$ be an IFS, let $p$ be a probability vector, and let $\mu$ be the
corresponding self-similar measure.
Then
\[
\dim \mu =\min\Big(1,\frac{h(\Phi,p)}{\sum_{i\in\Lambda}p_i\log|\lambda_i|^{-1}}\Big).
\]
\end{Conjecture}

A special case of this conjecture, when $\Phi$ is assumed to have no exact overlaps
predicts that the dimension of $\mu$ equals the minimum of $1$ and the similarity
dimension \eqref{eq:sim_dim}.
A version of this for self-similar sets was first formulated by Simon \cite{Sim}.

This conjecture is wide open, but there are some recent results in some special cases.
The purpose of this survey is to discuss the role of entropy rates in these developments.
For broader discussions of self-similar sets and measures, we refer to the 
recent book \cite{BSS} and the survey \cite{Var:ICM}.
In what follows we discuss results in settings of increasing arithmetic complexity.
We begin with IFS's where all scaling parameters are algebraic.
Then
we consider the cases where all scaling parameters are the same possibly transcendental number
but the translation parameters are algebraic.
Finally, we discuss some cases, where there is transcendence in both the scaling and
translation parameters.

\subsection*{Acknowledgements}
I thank Ariel Rapaport, Lauritz Streck and the referee for their careful reading of an earlier
version of this survey.

\section{Algebraic contractions}\label{sc:algebraic-contractions}

About a decade ago, Hochman \cite{Hoc:1d} has made a major breakthrough in the dimension theory of self-similar
measures, which has inspired much further work.
He proved that Conjecture~\ref{cn} holds under what has become known as the exponential separation
condition.

\begin{Definition}\label{df:exp-separation}
We define the distance $\dist(\f_1,\f_2)$ of two similarities $\f_1,\f_2:\R\to\R$ to be $\infty$ if they have different
scaling factors and to be $|\f_1(0)-\f_2(0)|$ otherwise.
An IFS $\Phi$ is said to satisfy the exponential separation condition if there is a constant $c>0$
such that for infinitely many $k\in\Z_{>0},$ we have that
$(i_1,\ldots,i_k)\neq(j_1,\ldots,j_k)\in\Lambda^k$ implies
\[
\dist(\f_{i_1}\circ\cdots\circ\f_{i_k},\f_{j_1}\circ\cdots\circ\f_{j_k})>0.
\]
We say that $\Phi$ satisfies the weak exponential separation condition if the condition
$(i_1,\ldots,i_k)\neq(j_1,\ldots,j_k)\in\Lambda^k$ is replaced by
$\f_{i_1}\circ\cdots\circ\f_{i_k}\neq\f_{j_1}\circ\cdots\circ\f_{j_k}$ in the definition.
\end{Definition}

Observe that the exponential separation condition implies that the IFS $\Phi$ has no
exact overlaps, while the weak exponential separation condition does not.
In fact, it is a version of the condition that is suitable for studying IFS's with
exact overlaps.

The (weak) exponential separation condition is much less restrictive than classical separation conditions
such as the strong separation condition or the open set condition.
Indeed, while those conditions mean that the images of $\supp \mu$ are (almost) disjoint under the
maps in the IFS, the exponential separation condition allows these to overlap non-trivially.
Moreover, the exponential separation condition holds outside a $1$-codimensional set of potential
exceptions in any reasonable parametrized family of IFS's (see \cite{Hoc:1d}*{Theorem 1.8}
and \cite{Hoc:d}*{Theorem 1.10} for precise statements)
and also when all scaling and translation parameters of the maps in the IFS are algebraic numbers.

The following result of Hochman \cite{Hoc:1d}
is a major progress towards Conjecture~\ref{cn}.

\begin{Theorem}[Hochman]\label{th:hochman}
Let $\Phi$ be an IFS that satisfies the weak exponential separation condition,
let $p$ be a probability vector, and let $\mu$ be the
corresponding self-similar measure.
Then
\[
\dim \mu =\min\Big(1,\frac{h(\Phi,p)}{\sum_{i\in\Lambda}p_i\log|\lambda_i|^{-1}}\Big).
\]
\end{Theorem}

This result was stated under the exponential separation condition in \cite{Hoc:1d}, but
it follows easily from Hochman's work under the weaker condition, too.
The details of this are written in \cite{BV:overlaps}*{Theorem 3.5}.
The main result of \cite{BV:overlaps} gives an explicit formula for the dimension of
a natural family of self-similar measures with exact overlaps as an application
of Theorem~\ref{th:hochman}.

Theorem~\ref{th:hochman} implies in particular that Conjecture~\ref{cn} holds for IFS's if all
scaling and translation parameters are algebraic numbers.
Rapaport \cite{Rap:self-similar} improved this substantially by showing that it is enough to assume the condition
for the scaling parameters.

 \begin{Theorem}[Rapaport]\label{th:rapaport}
Let $\Phi$ be an IFS that contains no
exact overlaps such that the scaling parameters of all the maps are algebraic numbers.
Let $p$ be a probability vector, and let $\mu$ be the
corresponding self-similar measure.
Then
\[
\dim \mu =\min\Big(1,\frac{H(p)}{\sum_{i\in\Lambda}p_i\log|\lambda_i|^{-1}}\Big).
\]
\end{Theorem}

One may probably drop the no exact overlaps condition and replace $H(p)$ by the entropy rate
in the above theorem.
This stronger version is likely within the scope of the methods of \cite{Rap:self-similar},
however, this in not in the literature.

The proof of Theorem~\ref{th:rapaport} relies on Theorem~\ref{th:hochman} but it cannot be obtained as
a direct corollary.
Indeed, Baker \cite{Bak1}
(see also \cite{Che}, \cite{Bak2}, \cite{BK})
constructed IFS's satisfying the conditions of Theorem~\ref{th:rapaport}
such that the smallest distance of
distinct $k$-fold compositions go to $0$ arbitrarily fast as a function of $k$, so the
exponential separation condition may fail dramatically.
The key to the proof of Theorem~\ref{th:rapaport} is showing that there are very few pairs of
such $k$-fold compositions that are super-exponentially close.
Informally speaking, if we were to collapse
those $k$-fold compositions that are of super-exponentially small distance then this would
result in a reduction of entropy that becomes negligible as $k$ goes to infinity.
This means that while the minimal distance between some compositions of maps may decay
very fast, this happens so rarely that we can ignore it.

\section{Bernoulli convolutions}\label{sc:Bernoulli}

For a given parameter $\l\in(0,1)$, the Bernoulli convolution $\nu_\lambda$ is
the self-similar measure corresponding to the similarities $x\mapsto \l x\pm 1$ and
probability weights $1/2,1/2$.
Their study goes back to \cite{JW}, \cite{Erd1}, \cite{Erd2} predating the general theory.

Combining the results of Hochman \cite{Hoc:1d} and Varj\'u \cite{Var:Bernoulli},
we know that Conjecture~\ref{cn} holds
for this family of self-similar measures.
Denoting by $h(\lambda)$ the entropy rate of the IFS $\{x\mapsto \l x\pm 1\}$
with weights $1/2,1/2$, we have the following result.

\begin{thm}[Hochman, Varj\'u]\label{th:Bernoulli}
We have
\[
\dim\nu_\lambda=\min(1,h(\lambda)/\log\l^{-1}).
\]
\end{thm}

Wu \cite{Wu} has an exciting new approach for this problem.
He proved, among other results, that Bernoulli convolutions have dimension $1$
for all but countably many parameters in $[1/2,1)$.
While his argument gives no information about the exceptional parameters, his
approach has the advantage that it applies for all $1$ parameter families of IFS's
as long as a transversality condition is satisfied.
The proof is also considerably simpler relying only on Theorem~\ref{th:hochman}
and an ingenious argument using transversality to show that points where
the exponential separation condition fail (in a suitably strong sense)
repel each other. 

\subsection{Exceptional parameters in $[1/2,1)$}\label{sc:exceptional}
It is also natural to ask the question that given a parametrized family
of self-similar measures what the set in the parameter space is for which
the dimension of the self-similar measure is strictly less than that of the ambient space.
The entropy rate appears to be a very complicated function of the parameters,
so it is not expected that this question has an easily describable answer in general.
However, in the special case of Bernoulli convolutions, there may be a nice answer.

It is easy to see that $h(\lambda)\le \log 2$, which is the Shannon entropy
of the probability weights $(1/2,1/2)$.
So we have $\dim\nu_\l<1$ for all $\l\in(0,1/2)$.
In fact, this can be seen easily without using Theorem~\ref{th:Bernoulli}.
Indeed, it is not hard to show that $\supp\nu_\l$ is a Cantor set of dimension
less than $1$ in this parameter range.

By Theorem~\ref{th:Bernoulli}, we know that for $\l\in[1/2,1)$, we have $\dim\nu_\l<1$ if and
only if $h(\lambda)<\log\l^{-1}$.
This may only happen if the IFS has exact overlaps, which means that there are distinct
choices of the signs $(\o_1,\ldots,\o_n),(\wt\o_1,\ldots,\wt\o_n)\in\{-1,1\}^n$
for some $n$ such that
\[
\o_n\l^{n-1}+\ldots+\o_2\l+\o_1=\wt\o_n\l^{n-1}+\ldots+\wt\o_2\l+\wt\o_1,
\]
which is equivalent to $\l$ being the root of a polynomial with coefficients $-1$, $0$
and $1$.
In particular $\l$ is an algebraic unit for all such parameters.
(A number that is a root of a polynomial with integer coefficients whose
leading and constant terms are $\pm1$ is called an algebraic unit.)

Theorem~\ref{th:Bernoulli}, therefore, reduces our question to finding all algebraic units
$\l\in[1/2,1)$ with $h(\lambda)<\log\l^{-1}$.
The only known examples of such numbers are the inverses of so-called Pisot numbers.
An algebraic integer $\b>1$ is called Pisot if all its Galois conjugates have modulus
less than $1$.
Erd\H os \cite{Erd1}
showed that if $\lambda\in(1/2,1)$ is the inverse of a Pisot number, then
$\nu_\l$ is singular with respect to the Lebesgue measure by showing that its Fourier
transform does not vanish at infinity.
Later Garsia \cite{Gar:entropy}
built on Erd\H os's work to show that $h(\l)<\log\l^{-1}$ for these parameters.

\begin{Question}
Are there any $\l\in[1/2,1)$ such that $\dim\nu_\l<1$ and $\l^{-1}$ is
not a Pisot number?
\end{Question}

The Mahler measure of an algebraic number is an important quantity in Diophantine analysis, and
it turns out to be closely related to entropy rates.
Let $\l$ be an algebraic number with (not necessarily monic) minimal polynomial $P(x)\in\Z[x]$.
Let $P(x)=a_n(x-\a_1)\cdots(x-\a_n)$ be the factorization of $P$.
Then the Mahler measure of $\l$ is defined as
\[
M(\l)=|a_n|\prod_{j=1}^n\max(1,|\a_j|).
\]

If $\l$ is an algebraic unit, then the leading coefficient of the minimal polynomial is $1$,
and the product of the Galois conjugates has modulus $1$.
Therefore,
\[
M(\l)=\prod_{j=1}^n\max(1,|\a_j|)=\prod_{j=1}^n\max(1,|\a_j|^{-1}).
\]
In particular, we always have $M(\l)\ge |\l|^{-1}$.
Observe that when $\l^{-1}$ is a Pisot number, then the only factor with $|\a_j^{-1}|>1$
is $|\l|^{-1}$, hence $M(\l)=\l^{-1}$ in this case.
There is another family of algebraic units with this property.
An algebraic integer $\beta>1$ is called a Salem number if it has at least one Galois
conjugate on the unit circle and all Galois conjugates off the circle are $\beta$ and $\b^{-1}$.
It is easy to see that $M(\l)=\l^{-1}$ if and only if $\l^{-1}$ is a Pisot or Salem number.

There is a simple upper bound for the entropy rate in terms of the Mahler measure.
It can be seen that there is a one-to-one correspondence between compositions of the
maps $\{x\mapsto\l x\pm1\}$ and the values of polynomials with coefficients $0$ and $1$
evaluated at $\l$.
The number of such values can be bounded by $M(\l)^{n+o(1)}$, where $n$ is the degree
of the polynomials.
This leads to the bound $h(\l)\le \log M(\l)$.

Generalizing the argument of Garsia \cite{Gar:entropy} for Pisot numbers,
Breuillard and Varj\'u \cite{BV:entropy} proved that
$h(\l)<\log M(\l)$ when $\l$ is an algebraic unit without conjugates  on the unit circle.
This does not imply the existence of new Bernoulli convolutions with $\dim\nu_\l<1$.
Indeed, when $\l^{-1}$ is a Salem number, there are conjugates on the unit circle.
If it is not a Salem or Pisot number then $\log M(\l)>\log(\l^{-1})$.
It appears to be a very delicate question to decide if $h(\l)=\log\l^{-1}$ for (any or all?)
Salem numbers, or more generally when does the equality $h(\l)=\log M(\l)$ hold.

On the other hand, among non-unit algebraic numbers, Garsia exhibited a family of examples
such that $h(\l)=M(\l)=\log 2$, which are known as Garsia numbers.
These exhibit very nice properties.
In particular, Garsia \cite{Gar:arithmetic}
showed that the corresponding Bernoulli convolution is not only
of dimension $1$ but also absolutely continuous with respect to the Lebesgue measure 
with bounded density, and Dai, Feng and Wang \cite{DFW}
showed that it has power Fourier decay.
Streck \cite{Str:maximal} gave a very satisfactory 
description of this phenomenon for more general self-similar
measures.

It is also possible to bound the entropy rate from below in terms of the Mahler measure.
Breuillard and Varj\'u \cite{BV:entropy} proved
\begin{equation}\label{eq:entropy-lower-bound}
h(\lambda)\ge 0.44\log \min(2, M(\l)).
\end{equation}
Note that $h(\lambda)\le \log 2$, so the lower bound matches the upper bound up to a multiplicative
constant.
This implies that at least in the interval $[2^{-0.44},1)$, an exceptional parameter
(that is one for which $\dim \nu_\l<1$) must be close to being  Pisot or Salem 
in the sense that $M(\lambda)\le \l^{-1/0.44}$.
A famous open problem known as Lehmer's conjecture asserts that $M(\lambda)>1+c$ for a
positive absolute constant $c$ whenever, $\lambda$ is an algebraic number that is not
a root of unity.
If this holds, then there is an interval $(x,1)$ with some $x<1$ that does not contain
any exceptional parameters.
Unconditionally, this is an open problem.

Given a specific algebraic unit $\l$, there is an algorithm
\cite{AFKP} for approximating $h(\l)$
with arbitrary precision, see also \cite{HKPS} and \cite{FF}
for numerical estimates for entropy rates.
Therefore, for any specific number such that $h(\lambda)\neq \lambda^{-1}$, it is a
finite computation to decide if the corresponding Bernoulli convolution is of dimension $1$.
It would not be unreasonable to expect that $h(\lambda)=\lambda^{-1}$ may occur
only for Salem numbers.

\subsection{Continuity properties of entropy rates}\label{sc:continuity}

It is well known (see \cite{HS}*{Theorem 1.8} or \cite{BV:dimension}*{Lemma 16}
for details) that the function $\lambda\mapsto \dim\nu_\l$
is lower semi-continuous on $(0,1)$.
This and Theorem~\ref{th:Bernoulli} implies that $\lambda\mapsto \min(\log \lambda^{-1},h(\lambda))$ is also
lower semi-continuous on $(0,1)$.
(Note that $h(\l)=\log 2$ for $\l\in(0,1/2]$.)
This provides some motivation to ask the following question.

\begin{Question}\label{qu:continuity}
Is it true that the map $\l\mapsto h(\l)$ is lower semi-continuous on $(0,1)$?
\end{Question}

Admittedly, this is highly speculative given that Theorem~\ref{th:Bernoulli} or its proof says nothing
about the values of $h$ larger than $\log\lambda^{-1}$, which is the case outside a tiny set.

This question is closely related to some well known but poorly understood questions
about Mahler measure.
Given some $n\in\Z_{\ge 2}$, denote by $h(\l,n)$ the entropy rate of the IFS
\[
\{x\mapsto \l x, x\mapsto \l x+1,\ldots, x\mapsto \l x +(n-1)\}
\]
with the probability weights $(1/n,\ldots,1/n)$.
A simple change of coordinates show that $h(\l)=h(\l,2)$.
If $\l\mapsto h(\l)$ was lower semi-continuous, then it would not be unreasonable
to expect that the functions $\l\mapsto h(\l,n)$ would also be lower semi-continuous
for each $n$.
It is not hard to show that $h(\l,2^n)$ monotone increases and converges to $\log M(\l)$
for all algebraic numbers $\l$.
(Indeed, this follows from the methods of \cite{BV:entropy}, for example,
but it could be proved in a simpler way.)
We could conclude hence that the map $\lambda\mapsto M(\l)$  was  lower semi-continuous
on the real line if we defined $M(\l)=\infty$ for transcendental $\l$.
If true, this would be a very beautiful result.
A special case of it is that the set
\begin{align*}
\{\l\in(0,1):M(\l)\le \log \l^{-1}\}&=\{\l\in(0,1):M(\l)= \log \l^{-1}\}\\
&=\{\l\in(0,1):\l^{-1}\text{ is Pisot or Salem}\}
\end{align*}
is closed.
This is a well known conjecture of Boyd \cite{Boy}*{Section 4}.

\subsection{A word about the proof of Theorem~\ref{th:Bernoulli}}\label{sc:proof}

Theorem~\ref{th:Bernoulli} for algebraic parameters is due to Hochman \cite{Hoc:1d},
see also \cite{BV:entropy}*{Section 3.4}.

Entropy rates play an important role in the proof even in the transcendental
case.
We outline this here briefly.
The proof is given in the short paper \cite{Var:Bernoulli},
but it builds on several major
previous developments.
One of these ingredients is a result of Breuillard and Varj\'u \cite{BV:dimension}
that if there is a putative transcendental $\lambda\in(1/2,1)$
with $\dim\nu_\l<1$, then for infinitely many $n$'s
there is an algebraic unit $\eta$ of degree at most $n$ such that
$\eta$ approximates $\l$ with very small error:
\[
|\eta-\l|<\exp(-n^2),
\]
and $h(\eta)<\log 2-c$ with a constant $c>0$ that depends only on $\l$
(but not on $n$).
This was proved by an argument that builds on Hochman's work \cite{Hoc:1d}.

The next step is to use the information $h(\eta)<\log 2-c$ to prove $M(\eta)< C$
with a constant $C$ that depends again only on $\l$.
This can be proved by an argument similar to the proof of \eqref{eq:entropy-lower-bound}
in \cite{BV:entropy}.

Next, a result in Diophantine analysis is used which implies that if
$m$ is sufficiently larger than $n$ ($>n\log n $ or $>n(\log n)^2$ is
enough depending on which argument is used), then a neighbourhood of $\eta$
that is exponentially small in $m$ does not contain another algebraic parameter
of degree at most $m$ for which the IFS has exact overlaps.
This can be done in several ways.
One argument is given in \cite{Var:Bernoulli}.
A simpler approach that also generalises easier is due to Dimitrov and appeared
in \cite{RV}*{Lemma 4.6, Appendix A}.

The conclusion of the previous step implies that the IFS satisfies the exponential separation
condition for a suitably chosen $m$ in the role of $k$ in Definition~\ref{df:exp-separation}.
Starting with different $n$'s we can find an infinite sequence of $m$'s with this property,
hence the IFS for this $\l$ satisfies the exponential separation property.
This contradicts to Theorem~\ref{th:hochman} and the assumption $\dim\nu_\l<1$.

\subsection{Generalizations}\label{sc:generalizations}

The argument can be generalized to other IFS's.
This was first done in \cite{RV}*{Appendix A} for an IFS of the form
\[
\{x\mapsto \l x+\tau_1,\ldots,x\mapsto \l x+ \tau_k\}
\]
where $\l\in(0,1)$ and $\tau_1,\ldots,\tau_k$ are rational numbers.
This was further extended in \cite{FF:algebraic} to the case where $\tau_1,\ldots,\tau_k$
are only assumed to be algebraic.

Furthermore, Rapaport \cite{Rap:d} and Rapaport and Ren \cite{RR} built on these techniques
to compute the dimension of the self-affine measures corresponding to an IFS consisting of maps
of the form
\begin{equation}\label{eq:self-affine}
\f_i:(x_1,\ldots,x_d)\mapsto (\l_1 x_1,\ldots,\l_d x_d) + a_i
\end{equation}
where $a_i$ are rational vectors assuming only that the coordinate projections
of the IFS have no exact overlaps.

\section{Two transcendental parameters}\label{sc:two-parameters}

Bernoulli convolutions form a family of IFS parametrized by a single
variable, the common scaling factor of the two maps.
The rest of the parameters are fixed algebraic numbers.
Rapaport and Varj\'u \cite{RV} have made progress towards extending the theory
to certain families parametrized by two variables.
For simplicity, we limit our discussion to the family of IFS
\begin{equation}\label{eq:2-parameters}
\Phi_{\l,\tau}=\{x\mapsto\l x\pm1\pm \tau\}.
\end{equation}
The signs may be chosen independently, so there are $4$ maps in the IFS,
and we will consider the probability weights $1/4,1/4,1/4,1/4$.
The paper \cite{RV} is more general, in particular, it also treats the
case of IFS's that consist of $3$ arbitrary similarities with a common
scaling factor.
We refer to the original paper for the details.

Let $n\in\Z_{\ge 1}$, and let 
\[
\f_1(\l,\tau),\ldots,\f_n(\l,\tau),\psi_1(\l,\tau),\ldots,\psi_n(\l,\tau)\in\Phi_{\l,\tau}
\]
be a choice of $2n$ elements from the family of IFS.
Expanding the equation
\begin{equation}\label{eq:curve}
\f_1(\l,\tau)\circ\ldots\circ\f_n(\l,\tau)=\psi_1(\l,\tau)\circ\ldots\circ\psi_n(\l,\tau)
\end{equation}
we get a polynomial equation in the the two variables $\l,\tau$.
In the $\tau$ variable, the equation is, in fact, linear.
The solution of this equation (assuming it is non-trivial) is a finite union of curves
in the $\l,\tau$ parameter space, for each of whose points, the IFS has exact overlaps.

This is a crucial difference compared to one parameter families like Bernoulli convolutions.
In that case, the parameters for which exact overlaps occur up to a certain length of compositions,
is a finite set.
The separation between these points plays a decisive role in the argument.
In the case of two parameters, exact overlaps occur along curves, which may intersect,
and hence there is no minimal
separation between parameter points with exact overlaps.

The case of IFS's with algebraic scaling and arbitrary translation parameters,
which are treated by Rapaport in \cite{Rap:self-similar}
as we mentioned in Section~\ref{sc:algebraic-contractions}
is also a multi-parameter family, and the same difficulty arises.
However, in that case, all the equations are linear, which makes this
case easier to handle than the family \eqref{eq:2-parameters}.

The case of the family of self-affine IFS's \eqref{eq:self-affine},
which we mentioned in Section~\ref{sc:generalizations}, also involves
multiple parameters.
However, in this case the number of parameters matches the dimension
of the ambient space, and hence an exact overlap is described
by a system of equations where the number of equations matches the
number of the parameters.
Therefore, exact overlaps occur in discrete points.
In this case, the challenge lies in the non-conformal nature of the IFS. 

For the purpose of the arguments in the proof of Theorem~\ref{th:Bernoulli},
we may ignore
parameters with exact overlaps if the entropy rate is not significantly smaller
than $\log\l^{-1}$.
In the approach of \cite{RV} to extend these argument to treat the
family \eqref{eq:2-parameters},
a key step is showing that those parameter points where the entropy rate of the
IFS drops significantly below $\log\l^{-1}$ are discrete.

To explain this properly, we need to extend the notion of entropy rates in such a way
that allows
to quantify the amount of exact overlaps that occur along a subset of the parameter
space simultaneously.

\begin{Definition}[Entropy rate for subsets in parameter space]\label{df:entropy-rate-set}
Let $\Phi_{\l,\tau}$ be the parametrized family of IFS's defined above.
Let $A$ be a set of the $(\l,\tau)$ parameters.
Let $X_{1,\l,\tau},X_{2,\l,\tau},\ldots$
be a sequence of independent uniformly distributed random elements of
$\Phi_{\l,\tau}$, and let
$Y_{n,\l,\tau}=X_{1,\l,\tau}\circ\ldots\circ X_{n,\l,\tau}$.
The compositions are understood as compositions of similarity transformations,
but we also think about $Y_n$ as a random map that takes the parameter $\l,\tau$ 
to a similarity transformation of $\R$.
We write $Y_{n}|_{A}$ for the random map $A\to\Sim(\R)$ that is obtained by
restricting $Y_n$ to $A$.
The entropy rate of the parameter set $A$ is defined as
\[
h(A)=\lim_{N\to\infty}\frac{1}{N}H(Y_N|_A),
\]
where $H(\cdot)$ denotes Shannon entropy.
\end{Definition}

Here again, the sequence in the definition is subadditive, and hence its limit exists
and is equal to the infimum.

It is easy to see that $h(A)$ is equal to $\log 4$ unless $A$ is a subset of the solution
set of an equation of the form \eqref{eq:curve} in parameter space.
We will refer to such solution sets as ``curves".

One of the main results of the paper \cite{RV} is that each parameter point
$(\l,\tau)$ such that the IFS $\Phi_{\l,\tau}$ has no exact overlaps, has a neighbourhood $U_\e$
for all $\e>0$ that does not intersect any curves $\gamma$ with entropy rate
$h(\gamma)<\log \l^{-1}-\e$.
Consider a point $(\eta,\s)\in U_\e$ with $h(\eta,\s)<\log\l^{-1}-\e$, and
consider all the equations of the form \eqref{eq:curve} that come from the exact overlaps
at $(\eta,\s)$.
Each of those equations vanishes at a curve that passes through $(\eta,\s)$
and whose entropy rate is larger than 
$\log\l^{-1}-\e>h(\eta,\s)$.
Therefore, these equations cannot all vanish along the same curve, and $(\eta,\s)$
must be the solution of a system of two independent equations.
This implies in particular that the set of these points is finite when we are
considering exact overlaps coming from compositions of a given length $n$.

This allows us to bypass a major difficulty for generalizing the argument described in
Section~\ref{sc:proof} for proving Conjecture~\ref{cn} in the setting of the IFS's $\Phi_{\l,\tau}$ defined above.
All of this programme was carried out in the paper \cite{RV}, except that an appropriate
generalization of the result in \cite{BV:entropy} that bounds the Mahler
measure of parameters whose entropy rate is bounded away from $\log 2$,
which we mentioned in Section~\ref{sc:proof}.
This generalisation is still outstanding.
More precisely, it was proved in \cite{RV} that Conjecture~\ref{cn}
holds for the family $\Phi_{\l,\tau}$ when the IFS has no exact overlaps
if the answer to the following question is affirmative.

\begin{Question}\label{qu:Rapaport-Varju}
Is it true that for any $\e>0$, there is $\d>0$ and $C<\infty$
such that for all parameter points $(\l,\tau)$ with
\[
h(\l,\tau)\le \min(\log 4,\log \l^{-1})-\e
\]
and $\l\in(\e,1-\e)$, we have $M(\l)<C$, or $(\l,\tau)$ is contained in a curve $\gamma$
with $h(\gamma)\le\min(\log 4,\log \l^{-1})-\d$?
\end{Question}

Since $h(\l,\tau)\le h(\gamma)$ for any curve $\gamma$ that passes through $(\l,\tau)$,
the assumption that $h(\l,\tau)\le \min(\log 4,\log \l^{-1})-\e$ alone does not even imply
that $\l$ is algebraic let alone a bound on its Mahler measure.
For this reason, it is important to allow the second alternative in the conclusion.

Question~\ref{qu:Rapaport-Varju} remains wide open.
A substantial progress was made on it by Streck in his PhD thesis \cite{Str:PhD}.

\begin{thm}[Streck]\label{th:Streck}
For every $\e>0$, there is a constant $C$ such that
if $h(\l,\tau)\le (3/2)\log 2-\e$ for some parameter point $(\l,\tau)$
then $M(\l)<C$ or $\tau=0$.
\end{thm}

The set of parameter points $(\lambda,0)$ for $\l\in(-1,1)$ form a curve whose entropy rate
is $\log 2$.
It turns out that this is the curve in parameter space of the
lowest entropy rate and there
is no other with entropy rate less than $(3/2)\log 2$.
Indeed, this is a consequence of the theorem.
There is a countable collection of curves that are of the form $\tau=\l^m$ for $m\in\Z$
that all have entropy rates equal to $(3/2)\log 2$.
Therefore, any improvement of Theorem~\ref{th:Streck} would need to account for these curves.

Theorem~\ref{th:Streck} also makes progress towards Conjecture~\ref{cn}.
Together with the results of \cite{RV} it implies the following.

\begin{thm}\label{th:Streck2}
Suppose $(\l,\tau)$ is such that the IFS $\Phi_{\l,\tau}$ defined in \eqref{eq:2-parameters}
has no exact overlaps and $\l>2^{-3/2}$.
Then the self-similar measure corresponding to $\Phi_{\l,\tau}$ and uniform
probability weights is of dimension $1$.
\end{thm}

The self-similar measures that occur in this result are convolutions of the
Bernoulli convolution $\nu_\l$ with a scaled copy of itself.
For this reason, the result already follows from Theorem~\ref{th:Bernoulli} if $\l\ge1/2$,
but it is completely new in the range $\l\in(2^{-3/2},1/2)$.

\bibliography{bib}

\end{document}